\documentclass[12pt]{article}
\usepackage{amsmath}
\usepackage{graphicx,psfrag,epsf}
\usepackage{enumerate}
\usepackage{natbib}
\usepackage{url} % not crucial - just used below for the URL 

\newcommand{\blind}{0}

\begin{document}

\def\spacingset#1{\renewcommand{\baselinestretch}%
	{#1}\small\normalsize} \spacingset{1}

\if0\blind
{
	\title{\bf Why bother with Bayesian $t$-tests?}
	\author{Fintan Costello \\ School of Computer Science and Informatics,\\ University College Dublin\\
		and \\
		Paul Watts \\ Department of Theoretical Physics,\\National University of
		Ireland  Maynooth\\ }
	\maketitle
} \fi

\if1\blind
{
	\bigskip
	\bigskip
	\bigskip
	\begin{center}
		{\LARGE\bf Why bother with Bayesian $t$-tests?}
	\end{center}
	\medskip
} \fi

\bigskip
\begin{abstract}
		Given the well-known and fundamental problems with hypothesis testing via classical (point-form) significance tests, there has been a general move to alternative approaches, often focused on the Bayesian $t$-test.  We show that the Bayesian $t$-test approach does not address the observed problems with classical significance testing, that Bayesian and classical $t$-tests are mathematically equivalent and linearly related in order of magnitude (so that the Bayesian $t$-test providing no further information beyond that given by point-form significance tests), and that Bayesian $t$-tests are subject to serious risks of misinterpretation, in some cases more problematic than seen for classical tests (with, for example, a negative sample mean in an experiment giving strong Bayesian $t$-test evidence in favour of a positive population mean).  We do not suggest a return to the classical, point-form significance approach to hypothesis testing.   Instead we argue for an alternative distributional approach to significance testing, which addresses the observed problems with classical hypothesis testing and provides a natural link between the Bayesian and frequentist approaches. 
\end{abstract}

\noindent%
\vfill

\newpage
\spacingset{1.45} % DON'T change the spacing!
\section{Introduction}

It is clear that classical or point-form significance testing has serious problems: many statistically significant experimental results fail to occur reliably in replications \citep[e.g.][]{camerer2018evaluating,open2015estimating,klein2018many,klein2014investigating}, the chance of getting a statistically significant $p$-value increases with sample size, irrespective of the presence or absence of a true effect \citep{thompson1998praise} and point-form null hypotheses are always false (and to quote Cohen, 2016: ``if the null hypothesis is always false, what's the big deal about rejecting it?'').     In an attempt to address these  problems various researchers have argued for a move to Bayesian hypothesis testing approaches, with a particular focus on generalisations of Jeffrey's Bayesian $t$-test, which involves a Bayes Factor comparison with a nested, point-form null hypothesis \citep{jeffreys1948theory,gonen2005bayesian,fox2006two,rouder2009Bayesian,wang2016simple,schmalz2021bayes}.   In this paper we show that this Bayesian $t$-test approach does not, in fact, address any of these problems with classical null hypothesis testing.  Instead, the Bayesian $t$-test involves comparison against a point-form null which we know is always false (so what's the big deal about getting evidence against it?);  the form of the Bayesian $t$-test means that  probability of getting Bayesian evidence against the null increases with sample size, irrespective of the presence or absence of a true effect;  and the Bayesian $t$-test gives results which are simply a linear transformation of those obtained in classical significance tests and so are necessarily subject to the same problems of replication and reliability as seen in classical tests.  We also show that the Bayesian $t$-test is subject to serious risks of misinterpretation, arguably more problematic than those seen for classical tests.   We demonstrate these points in detail below, beginning with a  derivation of the general Bayesian $t$-test and next showing that these problems all hold with this general form (and so hold for all specific instantiations).  We the argue that researchers should move to hypothesis testing relative to distributional rather than point-form nulls, and show that the distributional approach does not suffer from any of these issues and further, naturally relates the Bayesian and frequentist approaches to hypothesis testing.  We conclude by briefly giving our view of the roles that Bayesian and frequentist statistical tools can play in scientific research.
   
   \section{The Bayesian $t$-test}
  
    The Bayesian t-test is a specific type of Bayes Factor test,  originally developed by Jeffreys (1948) and extended in various different ways by a range of other researchers. The Bayes Factor 
    \[ BF_{10} = \frac{p(y|H_1)}{p(y|H_0)}\]
    is the ratio of the likelihood of  observed data $y$ under hypothesis $H_1$ to its likelihood under $H_0$: the higher the value of $BF_{10}$, the more Bayesian evidence our data $y$ gives in favour of $H_1$ and against $H_0$.   More specifically, hypotheses or models $H_0$ and $H_1$ are taken to be probability distributions with some parameter values fixed and some parameters following specified prior distributions, and $p(y|H)$ is estimated in terms of the density of the $H$ distribution integrated over those priors (the marginal likelihood).        A Bayesian $t$-test is a Bayes Factor test where the observed data takes the form of a $t$ statistic, and where the hypothesis $H_0$ is the null hypothesis used in the classical $t$-test.  Various forms of Bayesian $t$-test have been proposed in the literature: while each uses a different form of alternative hypothesis $H_1$,  all have the same statistic $t$ and the same null hypothesis $H_0$.  
    
    Here we present Bayesian $t$-test in the context of a one-sample test, but with a generic structure which covers all possible forms of alternative hypothesis $H_1$.  We assume an experiment involving $N$ measurements or observations $X$ with observed sample mean $\overline{X}$, degrees of freedom $\nu$ and sample variance
    \[ S^2 = \frac{1}{\nu} \sum (X-\overline{X})^2 \]  
    We assume that observations $X$ follow a normal distribution
    \[X \sim \mathcal{N}(\mu,\sigma^2)\]
    for unknown parameters $\mu$ and $\sigma$, which means that
    \[ \frac{\nu S^2}{\sigma^2} \]
    follows a $\chi^2$ distribution with $\nu$ degrees of freedom. We let
    \[d=\frac{\overline{X}}{S^2}\] 
    represent the sample effect in this experiment and define the variable 
    \[ t= \frac{\overline{X}}{S}\sqrt{N} = d \sqrt{N} \]

     Our null hypothesis  is that $\mu=0$: that observations $X$ follow the normal distribution
    \[X \sim \mathcal{N}(0,\sigma^2)\]
     with unknown variance $\sigma^2$ (the null hypothesis in a classical $t$-test), so hypothesis $H_0$ is  
   \begin{equation}
   \label{eq:H_0_distribution}
    H_0: \overline{X} \sim \mathcal{N}(0,\sigma^2/N)
   \end{equation}
   and so the variable
   \[ \frac{\frac{\overline{X}}{\sqrt{\sigma^2/N}}}{\sqrt{\frac{\nu S^2}{\sigma^2}/\nu} } = t \]
   follows a  $T$ distribution with $\nu$ degrees of freedom.

We take the alternative hypothesis $H_1$ to have parameters $\sigma$,   $m$ and $\sigma_{m}$ such that observations $X$ follow the Normal distribution
    \[X \sim \mathcal{N}(\mu,\sigma^2)\]
and $\mu$ itself follows the Normal distribution
    \[\mu \sim \mathcal{N}(m, \sigma_{m}^2)\]
so that for given values of $\sigma$, $\sigma_{m}$ and $m$ we see that $\overline{X}$ follows the distribution 
 \begin{equation*}
 \begin{split}
 \overline{X} &\sim \int_0^{\infty} \mathcal{N}(\mu,\sigma^2/N) \mathcal{N}(\mu|m ,\sigma_{m}^2) \mathrm{d}\mu\\
 &=\mathcal{N}(m, \sigma^2/N+ \sigma_{m}^2) =  \mathcal{N}\left(m, (\sigma^2/N)\left[1+ \frac{\sigma_{m}^2}{\sigma^2}N\right]\right) 
 \end{split}
 \end{equation*}
Defining $\sigma_\delta = \sigma_{m}/\sigma$ (so that $\sigma_\delta^2$ is the variance of the effect size $\delta=\mu/\sigma$) this means that our hypothesis $H_1$ is
 \begin{equation} 
 \label{eq:H_1_distribution}
 \begin{split}
 H_1: \overline{X} & \sim  \mathcal{N}(m, (\sigma^2/N)\left[1+ \sigma_\delta^2 N\right]) 
 \end{split}
 \end{equation}
so that for fixed values of $m$ and $\sigma_\delta$ the variable
    \[ \frac{t}{\sqrt{1+\sigma_\delta^2 N}} =\frac{\frac{\overline{X} - m}{\sqrt{(\sigma^2/N)\left[1+ \sigma_\delta^2 N\right]}} + \frac{m}{\sqrt{(\sigma^2/N)\left[1+ \sigma_\delta^2 N\right]}} }{\sqrt{\frac{\nu S^2}{\sigma^2}/\nu} }  \]
 follows a non-central $T$ distribution with $\nu$ degrees of freedom and non-centrality parameter
 \[ \frac{m}{\sqrt{(\sigma^2/N)\left[1+ \sigma_\delta^2 N\right]}} = \frac{\delta}{\sqrt{1/N+ \sigma_\delta^2}}\]
 
 Letting $T_v$ represent the standard (central) $T$ distribution with $\nu$ degrees of freedom and $T_\nu(\theta)$ represent the non-central $T$ distribution  $\nu$ degrees of freedom and non-centrality parameter $\theta$, we can thus express our two hypotheses $H_0$ and $H_1$ in equivalent forms as
 \begin{equation} 
 \label{eq:H_0_alt}
 H_0 : t \sim T_\nu
 	\end{equation}
 and
  \begin{equation} 
 \label{eq:H_1_alt} H_1 :  \frac{t}{\sqrt{1+\sigma_\delta^2 N}} \sim T_\nu\left(\frac{\delta}{\sqrt{1/N+ \sigma_\delta^2}}\right)
 \end{equation}

In the Bayesian approach the likelihoods $p(t|H_0)$ and $p(t|H_1)$ are taken to be equal to the density of these these distributions at the values given in $H_0$ and $H_1$.  Taking $f_\nu(t)$ to be the density of the standard $T$ distribution at $t$ and $f_\nu(t;\theta)$ to be the density of the non-central $T$ with parameter $\theta$ the density of $H_0$ at $t$ is 
\[ p(t|H_0) = f_\nu(t) \]
and the density of $H_1$ at
\[\frac{t}{\sqrt{1+\sigma_\delta^2 N}} \]
is 
\[ p(t|H_1) = \frac{1}{\sqrt{1+\sigma_\delta^2 N}} f_\nu\left(\frac{t}{\sqrt{1+\sigma_\delta^2 N}};\frac{\delta}{\sqrt{1/N+ \sigma_\delta^2}}\right)  \]
giving
\begin{equation}
\label{eq:gronau}
BF_{10} = \frac{\frac{1}{\sqrt{1+\sigma_\delta^2 N}} f_\nu\left(\frac{t}{\sqrt{1+\sigma_\delta^2 N}};\frac{\delta}{\sqrt{1/N+ \sigma_\delta^2}}\right)  }{f_\nu(t)}
\end{equation}
Equation \ref{eq:gronau} represents, for example, a one-sample instantiation of the Bayesian $t$-test of G{\"o}nen et al. \citep{gonen2005bayesian,gronau2019informed} and, taking $\delta =0$ and adding a $\chi^2(1)$ prior on $\sigma_\delta$, represents the one-sample JZS Bayesian $t$-test of \citet{rouder2009Bayesian}.  Other forms of the Bayesian $t$ test are produced by assuming different prior distributions for the parameters $\delta$ and $\delta_\sigma$ or by expanding the hypotheses in various ways (Jeffrey's original formulation, for example, involves splitting $H_1$ into three component hypotheses); all approaches, however, take some analog of these $H_0$ and $H_1$ distributions as their starting point, and so this presentation characterises the general Bayesian $t$-test. 

A core distinction between different forms of Bayesian $t$-test concerns the choice of value or prior for $\delta$ (and so for $m$, the mean for the distribution of $\mu$ in the alternative hypothesis $H_1$).  Default or local tests assume that $\delta$ is either equal to $0$ (a  delta distribution) or has a mean of $0$.  This implies that $m$ is also equal to or has a mean of $0$, and so in these tests both $H_0$ and $H_1$ assume the same mean for $\overline{X}$.  Informed or non-local tests, by contrast, assume that $\delta$ is equal to or distributed around some non-zero value chosen on the basis of prior knowledge in some way, so that $\overline{X}$ is assumed to have a different mean in $H_1$ than in $H_0$.  In the next section we discuss various problems of interpretation that arise with default or local tests.

\subsection{Problems of interpretation: default tests }

Two points are immediately evident from this general presentation of the Bayesian $t$-test.  First,  distribution $H_0$ is the null hypothesis distribution that underlies the classical or point-form  $t$-test (which also assumes that $\overline{X}$ is normally distributed around a mean of $0$ with variance $\sigma^2/N$).    The point-form null hypothesis, however, is always false: and so it is not clear what is to be gained by testing against it.  Second,  for default tests with $\delta=0$ (and so  $m=0$) hypotheses $H_0$ and $H_1$ differ only in the variance they assign to $\overline{X}$ (compare Equation \ref{eq:H_0_distribution} to Equation \ref{eq:H_1_distribution} with $m=0$).    Here there is a serious risk of misinterpretation, arising because researchers commonly take default Bayesian $t$-test results in favour of $H_1$ as giving evidence that the population mean differs from $0$ (that is, evidence of a significant effect).   This is clearly incorrect: in a default test $H_1$ assumes that the population mean \textit{is} $0$, and we cannot take evidence in favour of $H_1$ as evidence against this assumption.

 \subsection{Bayesian evidence but no real effect}
For any fixed value of $\sigma$, the variance of $\overline{X}$ in $H_0$ falls with rising $N$ to a limit of $0$, since that variance is $\sigma^2/N$.    This means that the probability of getting any value $\overline{X} \neq 0$ under $H_0$ similarly falls to $0$ with rising $N$; and so, for  any positive value $y$ and any sample effect $d \neq 0$, there exists some $N_0$ such that  $p(t|H_0) = p(d\sqrt{N}|H_0)  < y$  for all $N \geq N_0$.  
 
 For any fixed values of  $\sigma$, $m$ and $\sigma_m$, however, the variance of $\overline{X}$ in $H_1$ falls with rising $N$ to a limit of $\sigma^2_m$ (since that variance is $\sigma^2/N + \sigma^2_m$).  This means that for any sample effect $d \neq 0$ there will  thus exist some value $N_0$ such that $p(d\sqrt{N}|H_1)>p(d\sqrt{N}|H_0)$ holds for all  $N \geq N_0$.   Given this we see that 
\[ \lim_{N\rightarrow \infty} BF_{10} = \lim_{N\rightarrow \infty} \frac{ p(t|H_1)}{ p(t|H_0)} = \infty \]
necessarily holds for any value $\overline{X} \neq 0$ (and hence any $t \neq 0$) and  any required level of evidence in favour of the alternative hypothesis in a Bayesian $t$-test can be obtained with large enough sample size $N$, irrespective of the presence or absence of a true effect for both default and informed tests and  irrespective  of the choice of priors.

\subsection{Problems of interpretation: informed tests }

  Comparing Equations \ref{eq:H_0_distribution} and \ref{eq:H_1_distribution} we see that for informed tests with $m \neq 0$, $H_0$ and $H_1$ differ both in their assumed means and in their models of variance for $\overline{X}$.    This means that evidence against $H_0$ and in favour of $H_1$ may arise as a consequence of this difference in variance alone; again, this leads to a serious risk of interpretation, where researchers may  assume that evidence in favour of $H_1$ indicates that the population mean is closer to or more consistent with the alternative mean $m$ than the null mean $0$.  This is not the case. 
  
  It may be useful to give a concrete example of the problem.     Suppose we have a one-sample experiment with sample size $N=50$ (and so $\nu=49$) and that in our null hypothesis we assume $\delta=0$ (there is no effect) and our alternative hypothesis we assume $\delta=0.5$ (there is a medium-sized positive effect).  For our Bayesian analysis, we make the standard choice of simple unit-information prior for $\delta$ of $\sigma_\delta=1$.    Suppose we observe a medium-sized negative effect in our experiment of $d=-0.5$ (so that $t=-0.5\sqrt{50}$).  Then applying Equation \ref{eq:gronau}  we have a Bayesian $t$-test comparing $H_1$ to $H_0$ of
    \[BF_{10} = \frac{\frac{1}{\sqrt{1+\sigma_\delta^2 N}} f_\nu\left(\frac{t}{\sqrt{1+\sigma_\delta^2 N}};\frac{\delta}{\sqrt{1/N+ \sigma_\delta^2}}\right)  }{f_\nu(t)} = \frac{\frac{1}{\sqrt{1+50}} f_{49}\left(\frac{-0.5 \sqrt{50}}{\sqrt{1+50}}, \frac{0.5 }{\sqrt{1+1/50}} \right)}{f_{49}\left(-0.5\sqrt{50}\right)} \approx 25 \]
This is strong Bayesian evidence in favour of the alternative hypothesis $H_1$, and if we mistakenly assume that $BF_{10}$ gives evidence about the hypothesised effect $\delta$ (as opposed to the effect-plus-variance model $H_1$), we will be led to the nonsensical conclusion that observing a medium-sized \textit{negative} effect in our experiment gives us strong Bayesian evidence in favour of a medium-sized \textit{positive} effect. 

Note that we pick on this one-sample instantiation of the \citet{gonen2005bayesian}  $t$-test here only because of its clarity and simplicity of presentation: the general problem  (of negative results giving apparently strong evidence in favour of a positive hypothesis) applies for all informed or non-local Bayesian $t$ tests, and arises, as before, because the variance of $H_0$ falls to $0$ with rising $N$ while the variance of $H_1$ does not.

\subsection{Default Bayesian $t$-tests and classical $t$-tests are equivalent}

Our last point involves the relationship between the Bayesian and  classical $t$-tests.   It has long been observed that  Bayesian $t$-test evidence in favour of the alternative and classical point-form evidence against the null are essentially equivalent  \citep[to quote Jeffreys: ``As a matter of fact I have applied my significance tests to numerous applications that have also been worked out by Fisher’s, and have not yet found a disagreement in the actual decisions reached''; cited in][]{ly2016harold}; here we explain why this relationship holds.  

We first note that for large $\nu$ we have $T_\nu \approx \Phi$ (the standard $T$ distribution is well approximated by the standard Normal distribution) and that for the standard Normal distribution  the Mills ratio
\[M(x) =\frac{\Phi(-|x|)}{\phi(x)}\]
(the ratio of the  cumulative Normal function at $-|x|$ to the probability density at $x$) has the  well-known asymptotic approximation
\begin{equation}
\label{eq:mills_1}
M(x) \approx \frac{1}{|x|}
\end{equation}
which is relatively accurate for  $|x| > 3$ \citep[e.g.][pp. 43]{small2010expansions}.  This means that the classical $p$ value for a given $t$ is approximated by
\[ p = 2T_\nu(-|t|) \approx 2\Phi(-|t|) \approx  \frac{2\phi(t)}{|t|} \]
and substituting the expression for the standard Normal density 
\[ \phi(x) = \frac{1}{\sqrt{2\pi}} e^{-\frac{x^2}{2}} \]
and taking the log gives
\[ \log (1/p) \approx \frac{t^2}{2} + \log(|t|) +\log \left(\sqrt{\pi/2}\right)  \] 

For a default Bayesian $t$-test with $\delta=0$ we have
\[ BF_{10}=  \frac{\frac{1}{\sqrt{1+ \sigma_\delta^2 N}} f_\nu\left(\frac{t}{\sqrt{1+ \sigma_\delta^2 N}}\right)  }{f_\nu(t)} \approx \frac{\frac{1}{\sqrt{1+ \sigma_\delta^2N}} \phi\left(\frac{t}{\sqrt{1+ \sigma_\delta^2 N}}\right)  }{\phi(t)} \]
and again substituting and taking the log gives
\[ \log (BF_{10}) \approx   \frac{t^2}{2}\left(1  -\frac{1}{2(1+\sigma_\delta^2 N)}\right)-  \log \left(\sqrt{1+ \sigma_\delta^2 N}\right) \approx \frac{t^2}{2}- \log \left(\sqrt{1+ \sigma_\delta^2 N}\right) \]
and thus
 \[ \log (BF_{10}) \approx   \log (1/p) -\log(|t|) - \log \left(\sqrt{\pi(1+ \sigma_\delta^2 N)/2}\right)   \]
It is clear that large changes in the value of $t$  cause large changes in the value of $t^2/2$  but much smaller changes in $\log|t|$.  This means that if we have a set of experiments with approximately the same sample size $N$ (so that changes in the $\log(\sqrt{\pi (1+\sigma_\delta^2 N) /2})$ term across experiments are small) the we expect 
\[   \log BF_{10} \approx  \log 1/p + C \] 
to hold across a given set of experiments for some constant 
\[C = -\left\langle \log\left(|t| \sqrt{\pi (1+\sigma_\delta^2 N) /2}\right) \right\rangle  \]
where $\langle x \rangle$ indicates the average value of $x$ in those experiments.
This tells us that the Bayesian $t$-test $BF_{10}$ and the point-form significance $p$ are equivalent, at least in terms of order of magnitude.   Our main concern when considering statistical significance (or Bayesian evidence) is in the order of magnitude of our result rather than its exact value: in this context the Bayesian $BF_{10}$ and the classical $p$-value convey the same information, and the two tests are essentially the same.

\subsection{Testing the equivalence between $p$ and $BF$}

We tested this predicted relationship between Bayesian and classical $t$-tests using data from the first Many Labs replication project \citep{klein2014investigating}.  This  involved the replication of $16$ different experimental tasks investigating a variety of classic and contemporary psychological effects covering a range of different topics.  Each experiment was originally published in the cognitive or social psychology literature, and was replicated by researchers in around $36$ different sites. Of these $16$ tasks, $11$ involved independent $t$-tests: we downloaded the data on all experimental replications of these $11$ tasks ($396$ experiments in total)  and used the standard $R$ $t.test$ function \citep{Rmanual2021} to calculate the $t$-test $p$ and the $ttest.tstat$ function \citep[from the BayesFactor package,][]{bayesFactor2021} to calculate the Bayesian $t$-test $BF_{10}$ for each of these experiments. The R script for this analysis is  available  online (see Supplementary Materials). 

This particular form of Bayesian $t$-test is a default test assuming an alternative hypothesis $H_1$ with  $\delta=0$ and and with effect sizes distributed normally around $\delta$ with variance $\sigma_\delta^2$ which itself follows an inverse $\chi^2$ distribution with $1$ degree of freedom.  Under this prior $\sigma_\delta^2$ is distributed around $1$: this prior is therefore equivalent to, though slightly less informative than, the unit information prior $\sigma_\delta^2 = 1$ we used in our earlier example. 

Since the $t$-tests in this dataset were all independent two-sample tests with $N_1$ samples in one group and $N_2$ in the other, we took the effective sample size in each experiment to be
\[ N_{\textit{eff}} =  \frac{1}{1/N_1+ 1/N_2} \]
and calculated the value
 \[ - \log\left(|t| \sqrt{\pi (1+N_{\textit{eff}}) /2}\right)\]
  for each experiment in this dataset, and took $C$ to be the mean of these values, giving $C=-2.81$ for these experiments.  Our prediction is that $\log(BF_{10})$ and $\log(1/p)$ will have a linear relation in these experiments, with a slope of $1$ and an intercept of $C$. To test this prediction we took the $p$ and $BF_{10}$ values for each individual experiment and calculated the best-fitting  $\log(BF_{10})$ vs $\log(1/p)$ line relating these values.   The best-fitting line had a slope of $1.02$ and an intercept of $-2.81 \pm 0.015$ (see Figure $1$): the predicted value $C=-2.81$ fell within this (quite narrow) interval, confirming the predicted relationship.

 \begin{figure*}[h!] 
 	\begin{center}
 		\scalebox{0.8}{\includegraphics*[viewport= 0 0  750 350]{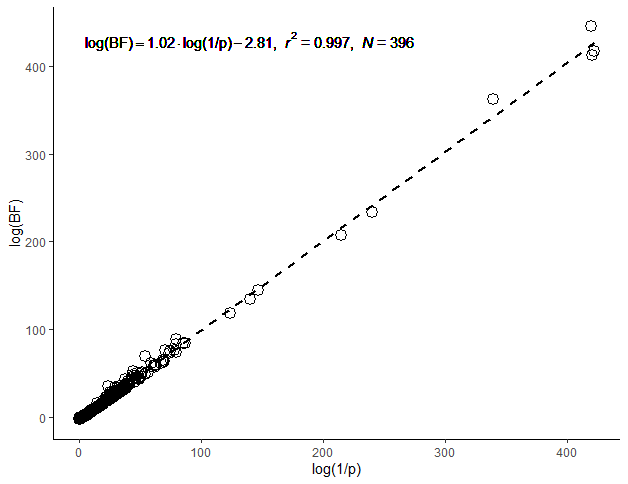}}
 	\end{center} 
 	\caption{ Scatterplot of $\log(BF_{10})$ vs $\log(1/p)$ for the $396$ $t$-test experiments in the Many Labs $1$ dataset, with the line of best fit.  The fit is extremely good, with the line accounting for more than $99\%$ of the variance in values; the slope is almost exactly the predicted value ($1.02$ vs $1$), and the intercept of $-2.81$ matches the predicted $C=-2.81$ value exactly. }
 \end{figure*}
 
 \section{Distributional null hypothesis testing}
 
 We've given a general characterisation of the Bayesian $t$-test and shown that this general form of the test, and so all specific instantiations, suffer from a series of problems: all compare an alternative $H_1$ against a null $H_0$ that we already know to be false; all give increasing evidence for $H_1$ irrespective of the presence or absence of any real effect; none give specific evidence about the population mean but instead give evidence about the variance of that mean; and (under a series of approximations) all are essentially equivalent to the classical $t$-test against a point-form null, providing no further information.  
 
 These problems arise from the use of the classical null hypothesis as $H_0$ in the Bayesian $t$-test approach, and from the fact that the two hypotheses $H_1$ and $H_0$ being compared differ in both their model of variance and (for informed tests) in their assumed mean. 
  Given these problems it seems unlikely that a move to Bayesian rather than classical hypothesis testing against the point-form null hypothesis will in any way address the problems with reliability and replication that we see in scientific research.  As an alternative, we suggest that researchers consider Fisherian evidential testing against a single null hypothesis,  but with a distributional rather than a point form null.  This is an approach where the statistical model is that observations $X$ follow the Normal distribution
\[X \sim \mathcal{N}(\mu,\sigma^2)\]
for unknown $\sigma^2$ and  where $\mu$ itself follows the Normal distribution
\[\mu \sim \mathcal{N}(m,\sigma_m^2)\]
  and where the null hypothesis is $m=0$.     We have recently proposed a distributional null hypothesis testing model following this approach which takes $\sigma_{m}^2$ to represent the variance in experimental means across replications of a given experiment.   The null is not always false in this model; evidence against the null in this model does not rise with sample size irrespective of the presence or absence of a real effect; and further, when the between-experiment variance of means is obtained from sample data, this model estimates the probability of replication of results in a way which reliably matches  observed rates of experimental replication  \citep[for a detailed presentation, see][]{costelloDistributional_t_submitted}.  

This distributional approach depends on a parameter $b = \sigma_{m}^2/\sigma^2$ representing the ratio of between-experiment variance in means to within-experiment variance in individual responses.  While this parameter is mathematically identical to the effect size variance $\sigma_\delta^2$ used in the derivation of the Bayesian $t$-test given above, it has a different meaning: where  $\sigma_\delta^2$ represents prior uncertainty about the effect (and so is subjective in nature), $b$ represents the relative variation in experimental means across different experiments and so is estimated from sample data (just as within-experiment variation $\sigma$ is estimated from sample data).  Further, where in a Bayesian $t$-test it is natural to choose an uninformative prior for $\sigma_\delta^2$, in the distributional approach the choice of value for $b$ represents a trade-off between Type $I$ and Type $II$ error: a high value for $b$ means high assumed between-experiment variance and so low Type $I$ error (but high type $II$ error), while a low assumed value for $b$ means low between-experiment variance and  so high Type $I$ error (but low type $II$ error)

This distributional null approach can also be applied to the comparison of null and alternative hypotheses; in this approach these two hypotheses are
\[ H_0:  \frac{t}{\sqrt{1+b N}} \sim T_\nu\left(\frac{t}{\sqrt{1+b N}}\right)  \]
and
\[ H_1:  \frac{t}{\sqrt{1+b N}} \sim T_\nu\left(\frac{t}{\sqrt{1+b N}};\frac{\delta}{\sqrt{1/N+ b}}\right)    
\]
and the Bayes Factor ratio for the alternative hypothesis $H_1$ against the null $H_0$ is the ratio of densities of these two distributions
\begin{equation*}
\begin{split}
 BF_{10} &= \frac{\frac{1}{\sqrt{ 1+ bN}}\, f_\nu\left(\frac{-|t|}{\sqrt{1+b N}};\frac{\delta}{\sqrt{1/N+ b}}\right)  }{\frac{1}{\sqrt{ 1+ bN}}\, f_\nu\left(\frac{-|t|}{\sqrt{1+b N}}\right)}   \\ 
 \end{split}
\end{equation*} 
and since both hypotheses $H_0$ and $H_1$ necessarily assume the same variance for $\overline{X}$ but different means, Bayesian evidence in favour of $H_1$ indicates that the observed data is more consistent with the mean in $H_1$ than the mean in $H_0$.   We can illustrate this using the same  one-sample experiment described earlier with sample size $N=50$ (and so $\nu=49$), a null hypothesis $\delta=0$ (there is no effect) an alternative hypothesis $\delta=0.5$ (there is a medium-sized positive effect) and assuming, purely for comparison purposes, a unit-information value of $b=1$.  This gives 
  \[BF_{10} = \frac{ f_\nu\left(\frac{-|t|}{\sqrt{1+b N}};\frac{\delta}{\sqrt{1/N+ b}}\right)  }{ f_\nu\left(\frac{-|t|}{\sqrt{1+b N}}\right)}  = \frac{ f_{49}\left(\frac{-0.5 \sqrt{50}}{\sqrt{1+50}}, \frac{0.5 }{\sqrt{1+1/50}} \right)}{ f_{49}\left(\frac{-0.5 \sqrt{50}}{\sqrt{1+50}}\right)} = 0.69 \]
and the test gives weak evidence in  favour of $H_0$, which is just as we would expect given that the observed result $d=-0.5$ is not strongly consistent with either $H_0$ or $H_1$, but is slightly more consistent with $H_0$.   

In the distributional approach the significance of a given result $t$ relative to $H_0$ is
\[ p_{sig}(t|H_0) = 2T_\nu\left(\frac{-|t|}{\sqrt{1+b N}}\right)    
\]
while its significance relative to an  alternative hypothesis of some effect size $\delta$ is
\[ p_{sig}(t|H_1) = 2T_\nu\left(\frac{-|t|}{\sqrt{1+b N}};\frac{\delta}{\sqrt{1/N+ b}}\right)    
\]
and there is a linear relationship between these measures of significance and the Bayes Factor measures of evidence.  
To see this we approximate both these $T$ distributions with corresponding Normal distributions 
(a rough approximation since it takes the Normal to approximate the non-central $T$) giving
\[ p_{sig}(t|H_0) \approx 2\Phi\left(\frac{-|d|}{\sqrt{1/N+b}}\right)    
\]
and
\[ p_{sig}(t|H_1) \approx 2\Phi\left(\frac{-|d-\delta|}{\sqrt{1/N+ b}}\right)    
\]
and  relating to Mills ratio  gives
\begin{equation*}
\begin{split}
BF_{10} & \approx \frac{ p_{sig}(t|H_1)}{  p_{sig}(t|H_0)}\,  \frac{M\left(\frac{-|d|}{\sqrt{1/N+b }}\right)}{ M\left(\frac{-|d-\delta|}{\sqrt{1/N+b}}\right) }  
\end{split}
\end{equation*} 

  We are primarily concerned here with results $d$ which are close to $0$ or to $\delta$ (giving evidence in favour of one hypothesis or the other), and so have the Mills ratio argument $x$ approaching $0$ for one or other hypothesis.  The approximation in  Equation \eqref{eq:mills_1} diverges as $x \rightarrow 0$, however; and so, since $\phi(0)=\sqrt{\pi/2}$, we use the modified approximation 
\[M(x) \approx \frac{1}{\sqrt{2/\pi}+|x|} \]
which is relatively close to $M(x)$ for $x < 3$ and which asymptotically approaches Equation \eqref{eq:mills_1} (and so $M(x)$) as $x \rightarrow \infty$.
Given this we see that the Bayes Factor and the ratio of distributional significance have the approximate relationship
\begin{equation*}
\begin{split}
BF_{10} &  \approx \frac{ p_{sig}(t|H_1)}{  p_{sig}(t|H_0)} \frac{\sqrt{2b/\pi}+|d-\delta| }{\sqrt{2b/\pi}+|d|  }  \\ 
\end{split}
\end{equation*}  
and the Bayes Factor measure of relative evidence for $H_1$ over $H_0$ given by result $t$ is, to a first approximation, simply a linear transformation of the distributional significance ratio for $t$ under $H_1$ and $H_0$. 

\section{Discussion}

Our focus so far has been on the relationship between Bayesian and frequentist approaches to hypothesis testing in a particularly simple situation: the $t$ test.  Here we briefly discuss the relationship between these two approaches more generally.  We take as our starting point a account of the Bayesian/frequentist distinction as given in a recent primer on Bayesian statistics: 
\begin{quote}
	\textit{The key difference between Bayesian and frequentist inference is that frequentists do not consider probability statements about the unknown parameters to be useful.	Instead, the unknown parameters are considered to be fixed; the likelihood is the conditional probability distribution	$p(y|\theta)$ of the data ($y)$, given fixed parameters ($\theta$). In Bayesian inference, unknown parameters are referred
	to as random variables in order to make probability statements	about them. The (observed) data are treated as	fixed, whereas the parameter values are varied; the likelihood is a function of $\theta$ for the fixed data $y$.} \begin{flushright} \citep[][p. 7]{van2021bayesian} \end{flushright}
\end{quote} 
We expand on this account by noting that in both frequentist and Bayesian approaches we have some theory  of the generative process producing data $y$.  This theory gives us two things: first, an overall statistical model $M$ with some set of independent parameters $\theta$ such that  $y$ is assumed to follow the distribution $y \sim M(\theta)$; and second, a list $H$ containing, for each parameter $\theta_i$, a particular selected value for that parameter (with some associated uncertainty or variance in that value).  The variances associated with parameter values in $H$ allows us to distinguish between fixed and free parameters in our theory.  A parameter value is fixed by theory, in this view, when our theory requires a specific value for $H_i$ so that any change to that value would necessarily require us to abandon the theory: the variance of a fixed parameter is thus necessarily $0$ in this theory.   If a parameter is not fixed it is free, and its value must be estimated from data in some way, so that any such estimated value, and any change in that value, remains consistent with our theory (and such that the current best estimate, and its variance, is given by $H_i$).

 Both frequentist and Bayesian approaches typically assume the overall model $M$ is fixed, and consider either testing or updating values of the parameter values $H$ (``the'' hypothesis).    Frequentist inference considers $p(y| H)$: the probability distribution for data $y$ conditional on $\theta = H$ (on the assumption that the parameters are as described in $H$).   Bayesian inference considers $p(H_{new}|y,H)$: the updated parameter descriptions $H_{new}$, conditional on  data $y$ and on the prior values $H$.  This common structure means that both forms of statistical inference fall within a single unified framework defined by $M$ and $H$: any Bayesian prior $H$ can be tested via the frequentist inference $p(y|H)$ and any frequentist hypothesis  $H$ about $\theta$ can be updated via the Bayesian inference  $p(H_{new}|y,H)$.  Indeed both forms of inference can be applied to the same data $y$, by asking  whether $y$ is consistent with $H$ and, if not, updating to produce a more consistent description $H_{new}$ and then asking whether $y$ is consistent with this new description  \citep[these are the prior and posterior predictive checks commonly recommended in standard Bayesian workflows, even though these checks involve a frequentist hypothesis test; see e.g.][]{gelman2020bayesian,schad2021toward}.   Note that for  parameter values $H_i$ with variance of $0$ (fixed by theory) this updating process will never cause any change in $H_i$.  This means that if $p(y|H_{new})$ is less than some significance criterion $\alpha$, we can conclude that data $y$ is inconsistent with our overall theory in some way: that updating to produce a set of parameter values consistent with the data would either require us to change some values that are fixed in that theory,  to abandon our prior estimates $H$ for some or all of those values (which by assumption were consistent with that theory) or to abandon our statistical model $M$.
 
  It is necessarily the case in this unified framework that any form of hypothesis testing (that is, any situation where we ask whether data $y$ is consistent with some $H$) will necessarily involve the frequentist inference  $p(y| H)$ in some way.  It should not be surprising, therefore, that the default Bayesian $t$ test and the classical $t$ test are essentially equivalent; the equivalence arises because both depend on inferences of the form $p(y| H)$.

  \bigskip
  \bigskip
  \begin{center}
  	{\large\bf SUPPLEMENTARY MATERIAL}
  \end{center}
  
  The R script used in this paper available  at \url{https://osf.io/qajvu}, and automatically downloads the Many Labs $1$ dataset, carries out the analysis, and generates  Figure $1$.

   \bibliographystyle{Chicago}   
   \bibliography{references}

    \end{document}